\documentclass[a4paper,11pt]{article}
\pdfoutput=1 

\usepackage{jheppub} 

\usepackage[T1]{fontenc} 

\title{\boldmath An Algorithm for Solving Solvable Polynomial Equations of Arbitrary Degree by Radicals}

\author[a,b]{Song Li,}

\affiliation[a]{CAS Key Laboratory of Theoretical Physics, Institute of Theoretical Physics, Chinese Academy of Sciences, Beijing 100190, China}
\affiliation[b]{School of Physics Sciences, University of Chinese Academy of Sciences,  Beijing 100049, China}

\emailAdd{lisong@itp.ac.cn}

\abstract{This work provides a method(an algorithm) for solving the solvable unary algebraic equation $f(x)=0$ ($f(x)\in\mathbb{Q}[x]$) of arbitrary degree and obtaining the exact radical roots. This method requires that we know the Galois group as the permutation group of the roots of $f(x)$  and the approximate roots with sufficient precision beforehand. Of course, the approximate roots are not necessary but can help reduce the quantity of computation. The algorithm complexity is approximately proportional to the 4th power of the size of the Galois group of $f(x)$. The whole algorithm doesn't need to deal with tremendous polynomials or reduce symmetric polynomials.}

\begin{document} 
    \maketitle
    \flushbottom

    \section{Introduction}
    \label{sec:intro}
    It is universally known that an irreducible polynomial equation with rational coefficients can be solved by radicals if and only if its Galois group is solvable. Polynomial equations with this property are said to be solvable. However, there is no efficient general algorithm for solving all solvable polynomial equations. Usually, we use specific methods to solve equations of specific degrees, such as solvable quintics~\cite{dummit1991solving} and solvable sextics~\cite{HAGEDORN2000704}. On the other hand, we have quite a few efficient methods for obtaining Galois groups of irreducible polynomials\cite{stauduhar1973determination,GEISSLER2000653,durov2006computation1,durov2006computation2,fieker_kluners_2014,hauenstein2018numerical}, and many, even ancient, methods for approximating the zeros of unary polynomial equations. This work provides a method for obtaining exact radical roots using numerical zeros and Galois groups. The numerical zeros are not indispensable for this algorithm, but they can reduce the quantity of calculation greatly.

    \section{The Method without Numerical Roots}
    \label{sec:PureMathMethod}

    Let $f(x)\in\mathbb{Q}[x]$ be an irreducible polynomial. We can convert $f(x)=0$ into an equation with integer coefficients:
    \begin{equation}
	    a_nx^n+a_{n-1}x^{n-1}+\cdots+a_1x+a_0=0.
    \end{equation}
    Multiplying both sides of the equation by $a_n^{n-1}$ and taking $y=a_n x$, we get a monic irreducible polynomial equation with integer coefficients. Therefore, we only consider the case where $f(x)$ is a monic irreducible polynomial with integer coefficients in this work. Let $G$ denote the Galois group of $f(x)$, which is the permutation group of the zeros of $f(x)$. We need $G$ to be solvable (otherwise out of our scope), so we can get the composition series of $G$: 
    \begin{equation}
	    G=G_m\rhd G_{m-1}\rhd G_{m-2}\rhd\cdots\rhd G_1\rhd G_0=\{e\},
	\end{equation}
	where $G_i/G_{i-1}$ is a cyclic group whose order is prime $p_i$, and the generator is $\sigma_iG_{i-1}$. Suppose the zeros of $f(x)$ are $x_1,x_2,\cdots,x_n$. Here we take $x_1$ as an example to show how to find its radical formula.
	
	$G_1$ is a cyclic group with prime order. It acts on $x_1$ to form a orbit whose length is either 1 or $p_1$. If $G_1$ does not change $x_1$, we can skip $G_1$ directly and take $\theta_{0}(x_1,x_2,\cdots,x_n)=x_1$. Otherwise, we take $\zeta_1$ as the primitive $p_1$th root of unity, and make the Lagrange resolvents:
	\begin{equation}
	    \begin{split}
	        (\zeta_1^0,x_{1}) &= x_{1}+\sigma_1 x_{1}+\sigma_1^2x_{1}+\cdots+\sigma_1^{p_1-1}x_{1}, \\
	        (\zeta_1^1,x_{1}) &= x_{1}+\zeta_1^1\sigma_1x_{1}+\zeta_1^2\sigma_1^2x_{1}+\cdots+\zeta_1^{p_1-1}\sigma_1^{p_1-1}x_{1}, \\
	        (\zeta_1^2,x_{1}) &= x_{1}+\zeta_1^2\sigma_1x_{1}+\zeta_1^4\sigma_1^2x_{1}+\cdots+\zeta_1^{2(p_1-1)}\sigma_1^{p_1-1}x_{1}, \\
	        &\,\,\,\vdots\\
	       (\zeta_1^{p_1-1},x_{1}) &= x_{1}+\zeta_1^{p_1-1}\sigma_1x_{1}+\zeta_1^{2(p_1-1)}\sigma_1^2x_{1}+\cdots+\zeta_1^{(p_1-1)^2}\sigma_1^{p_1-1}x_{1}.
	    \end{split}
	\end{equation}
	In this work, the Galois group is considered as the permutation group of the roots. Therefore it does not act on the constant $\zeta_1$. We have $\sigma_1(\zeta_1^i,x_1)=\zeta_1^{-i}(\zeta_1^i,x_1)$, so $(\zeta_1^i,x_1)^{p_1}$ is $G_1$ invariant (of course, $(\zeta_1^0,x_1)$ itself is $G_1$ invariant, no need for $p_1$ power). Expanding $(\zeta_1^i,x_1)^{p_1}$ into a polynomial of $\zeta_1$ and simplifying it only with the relation $\zeta_1^{p_1}=1$, we get
	\begin{equation}
	    \begin{split}
	        (\zeta_1^0,x_{1}) &=  x_{1}+\sigma_1x_{1}+\sigma_1^2x_{1}+\cdots+\sigma_1^{p_1-1}x_{1},\\
	        (\zeta_1^1,x_1)^{p_1} &= \theta_{0}+\theta_{1}\zeta_1+\theta_{2}\zeta_1^2+\cdots+\theta_{(p_1-1)}\zeta_1^{p_1-1}, \\
	        &\,\,\,\vdots\\
	        (\zeta_1^{p_1-1},x_1)^{p_1} &= \theta_{p_1(p_1-2)}+\theta_{(p_1(p_1-2)+1)}\zeta_1+\\
	        &\qquad\theta_{(p_1(p_1-2)+1)}\zeta_1^2+\cdots+\theta_{(p_1(p_1-1)-1)}\zeta_1^{p_1-1},
	    \end{split}
	\end{equation}
	where each $\theta_{j}$ is an integer coefficient polynomial of some $\sigma^kx_1$, which is essentially an integer coefficient polynomial of the roots of the equation $f(x)=0$. Therefore we can calculate the action of $G$ on each $\theta_{j}$. Since $(\zeta_1^i,x_1)^{p_1}$ is $G_1$ invariant, each $\theta_{j}$ is also $G_1$ invariant. In fact, most of these $\theta_{j}$ are the same. For example, take $\tilde{\zeta}_1=\zeta_1^k$, which is also a primitive $p_1$th root of unity, so $\tilde{\zeta}_1$ has no essential difference from $\zeta_1$, therefore
    \begin{equation}
	    (\zeta_1^k,x_1)^{p_1} = (\tilde{\zeta}_1,x_1)^{p_1}=\theta_{0}+\theta_{1}\tilde{\zeta}_1+\theta_{2}\tilde{\zeta}_1^2+\cdots+\theta_{(p_1-1)}\tilde{\zeta_1}^{p_1-1}.
	\end{equation}
	Hence, $\{\theta_{(k-1)p_1},\theta_{((k-1)p_1+1)},\cdots,\theta_{(kp_1-1)}\}$ and $\{\theta_{0},\theta_{1},\cdots,\theta_{(p_1-1)}\}$ differ by only one permutation. More precisely, if $t$ is the modular inverse of $k$ mod $p_1$, satisfy $kt\equiv1({\rm mod}\,p_1)$. In order to find out which $\theta_{m}$ the $\theta_{((k-1)p_1+j)}$ is equal to, we calculate $tj\equiv l({\rm mod}\,p_1)$, where $0\leq l\leq p_1-1$.Thus,
	\begin{equation}
	    \zeta_1^{j}=\left(\zeta^{kt}\right)^{j}=\tilde{\zeta}_1^{tj}=\tilde{\zeta}_1^{\,l}.
	\end{equation}
	Therefore we get $\theta_{((k-1)p_1+j)}=\theta_{l}$.
	
	Further, we have
	\begin{equation}
	    (\zeta_1^0,x_{1})^{p_1}=\theta_{0}+\theta_{1}+\cdots+\theta_{(p_1-1)}.
	\end{equation}
	Therefore, as long as $\theta_{0},\theta_{1},\cdots,\theta_{(p_1-1)}$ are calculated, we can  get the value of each $(\zeta_1^i,x_1)$ by using $\zeta_1$ and the $p_1$th root extraction, and then obtain $x_1$:
	\begin{equation}
	    x_1=\frac{(\zeta^0_1,x_1)+(\zeta^1_1,x_1)+\cdots+(\zeta_1^{p_1-1},x_1)}{p_1}.
	\end{equation}
	The other $\sigma_1^kx_1$ (they are also the roots of $f(x)=0$) can also be obtained by the Lagrange resolvents:
	\begin{equation}\label{eq:sigmax}
	    \sigma^kx_1=\frac{(\zeta^0_1,x_1)+\zeta_1^{-k}(\zeta^1_1,x_1)+\cdots+\zeta_1^{-k(p_1-1)}(\zeta_1^{p_1-1},x_1)}{p_1}.
	\end{equation}
	
	In addition, each $\theta_{j}$ is an integer coefficient polynomial of degree $p_1$ of $x_1,x_2,\cdots,x_n$, and there is no term of degree lower than $p_1$. 
	
    Now our goal has become to derive the exact radical formula for each $\theta_{j}$. We can completely repeat the above process: if $\theta_{j}$ remains unchanged under the action of $G_2$, we can directly define $\theta_{j,0}=\theta_{j}$; otherwise, we construct the corresponding Lagrange resolvents using the primitive $p_2$th root of unity... In general, this process will loop until all $\theta_{j_1,j_2,\cdots,j_m}$ are formed. For each $i\leq m$, $\theta_{j_1,j_2,\cdots,j_i}$ has the following properties:
    \begin{enumerate}
        \item $\theta_{j_1,j_2,\cdots,j_i}$ is an integer coefficient polynomial of degree $p_1p_2\cdots p_i$ of $x_1,x_2,\cdots,x_n$, and there is no term of degree lower than $p_1p_2\cdots p_i$. In particular, the degree of $\theta_{j_1,j_2,\cdots,j_m}$ is $|G|$; 
        \item $\theta_{j_1,j_2,\cdots,j_i}$ is the invariant of $G_i$;
        \item For each $i\leq m$, there can only be at most $p_1p_2\cdots p_i$ different $\theta_{j_1,j_2,\cdots,j_i}$;
        \item When we know the exact radical formula of each $\theta_{j_1,j_2,\cdots,j_i}$, we can get the radical formula of each $\theta_{j_1,j_2,\cdots,j_i}$ with the help of $\zeta_i$ and the $p_i$th root extraction. However, there are phase uncertainties due to the $p_i$th root extractions. I will solve this problem in the next section.
    \end{enumerate}
    
    The loop ends after finding the expression of $\theta_{j_1,j_2,\cdots,j_m}$, thus the final question is how to find the value of these $\theta_{j_1,j_2,\cdots,j_m}$. Since $\theta_{j_1,j_2,\cdots,j_m}$ is $G$ invariant, hence $\theta_{j_1,j_2,\cdots,j_m}\in\mathbb{Q}$. In fact, $\theta_{j_1,j_2,\cdots,j_m}$ all are integers. This is because $f(x)$ is a monic integer coefficient polynomial, so the elementary symmetric polynomials $\tau_1,\tau_2,\cdots,\tau_n$ of $x_1,x_2,\cdots,x_n$ take integer values. Let $P$ denote the set of representative elements of the coset cluster $S_n/G$ of $G$ in the $n$th permutation group $S_n$, and $\theta$ is any integer coefficient polynomial of $x_1,x_2,\cdots,x_n$, then
    \begin{equation}
	    F(x)=\prod_{\sigma\in P}(x-\sigma\theta)
	\end{equation}
	will be $S_n$ invariant. Thus, the coefficients of $F(x)$ are symmetric polynomials of $x_1,$ $x_2,$ $\cdots,$ $x_n$ with integer coefficients, which can be expressed as integer coefficient polynomials of $\tau_1,$ $\tau_2,$ $\cdots,$ $\tau_n$. When the values of $x_1,x_2,\cdots,x_n$ are substituted, $F(x)$ will become a monic integer coefficient polynomial, and its rational root must be an integer. $\theta$, as the rational root of $F(x)$, will be an integer. 
	
	If $F(x)$ is found, it will not be difficult to search its integer roots. Thus, through the whole algorithm, we can obtain the exact radical formula of the zeros of any solvable polynomial in principle. Unfortunately, $F(x)$ is a polynomial of degree $n!/|G|$, and its coefficients are polynomials of $x_1,x_2,\cdots,x_n$ which the degrees can be up to $n!$. Therefore, it is generally impossible to find out $F(x)$ exactly, and we need other method to get the value of each $\theta_{j_1,j_2,\cdots,j_m}$.

    \section{Reduce the Computing Burden by Numerical Roots}
    \label{sec:NumRootMethod}
    
    Suppose we have obtained sufficiently accurate values for all the roots of $f(x)=0$, and know the action of the group $G$ on them (in other words, know which value corresponds to $x_1$ and which value corresponds to $x_2$, and so on), then we substitute these values into the expression of $\theta_{j_1,j_2,\cdots,j_m}$, we will get a value that is very close to an integer. Find the nearest integer value for this value, then we get the exact value of $\theta_{j_1,j_2,\cdots,j_m}$. More importantly, with these numerical roots, we can determine the uncertain phase of $\sqrt[p_i]{(\zeta_i^k,\theta_{(i-1)j})^{p_i}}$. Therefore, the whole root-finding procedure no longer contains uncertainty. 
    
    If we have obtained the group $G$ and the approximate values of the roots, but we have not yet found the corresponding relationship between the subscripts of the roots and the approximate values. At this point, we can obtain the correspondence by virtue of the property that $G$ invariant integer coefficient polynomials have integer values. To this end, we take the representative element set $P$ of $S_n/G$, and construct a stabilizer $\theta$ of $G$. We randomly set a corresponding relationship at first, and then for each $\sigma\in P$, calculate whether $\sigma\theta$ is close to an integer value. If so, the corresponding relationship between the subscripts of the roots and the approximate values can be determined by this $\sigma$. However, this idea does not always work. Finding the subscripts of the roots is not the subject of this work, it will not be discussed further here. Let's go back to the topic. 
    
    The degree of each term in $\theta_{j_1,j_2,\cdots,j_m}$ is $|G|$. Taking each $\zeta_i$ and each $x_j$ equal to 1, we will find that the sum of the coefficients of $\theta_{j_1,j_2,\cdots,j_m}$ is less than $p_1^{p_1p_2\cdots p_m}p_2^{p_2p_3\cdots p_m}\cdots p_m^{p_m}$. Let $N$ denote the sum of the coefficients, we have 
    \begin{equation}
	    \theta_{mj}(x_1,x_2,\cdots,x_n)\sim Nx^{|G|}.
	\end{equation}
	If $x=x_0+\Delta x$, where $x_0$ denotes the exact root, $\Delta x$ represents the error between the numerical root and the exact root. Hence,
	\begin{equation}
	    Nx^{|G|}=N(x_0+\Delta x)^{|G|}=Nx_0^{|G|}\left(1+\frac{\Delta x}{x_0}\right)^{|G|}\approx Nx_0^{|G|}+N|G|x_0^{|G|-1}\Delta x.
	\end{equation}
    In order to make this value close to the integer we are looking for, we need $\left|N|G|x_0^{|G|-1}\Delta x\right|<1/2$, then 
    \begin{equation}\label{eq:accuate}
	    \left|\frac{\Delta x}{x_0}\right|<\frac{1}{2N|G||x_0|^{|G|}}.
	\end{equation}
	This formula gives us a bad signal. When the number of the prime factors of |G| is large, $N$ will be enormous, which will place very high requirements on the accuracy of the numerical root. However, even for a numerical root-finding method as simple as bisection method, the error decreases exponentially with the number of program loop steps. Therefore, the algorithm presented here still has the complexity advantage, even though $N$ will be relatively large.
	
	We can also reduce the requirement for numerical root precision by adding a part of the time or space complexity. To this end, we treat terms with different coefficients in $\theta_{j_1,j_2,\cdots,j_m}$ separately, and ignore the common coefficients temporarily.  This is because terms with the same coefficients form the $G$ invariant, which is still an integer coefficient polynomial after ignoring the integer common factors. This operation can reduce the value of $N$. Even more extreme, we take a term $x_1^{k_1}x_2^{k_2}\cdots x_n^{k_n}$, where $k_1+k_2+\cdots+k_n=|G|$ and $k_1\geq k_2,\cdots ,k_n$. And suppose $H\subset G$ is the stable subgroup of $x_1^{k_1}x_2^{k_2}\cdots x_n^{k_n}$, then
	\begin{equation}
	    \frac{1}{|H|}\sum_{\sigma\in G}\sigma\left(x_1^{k_1}x_2^{k_2}\cdots x_n^{k_n}\right)
	\end{equation}
	is an invariant polynomial of $G$ with coefficients of 1, the sum of its coefficients is $N=|G|/|H|\leq|G|$. By finding the integer values of all these polynomials, we can get the values of all $\theta_{j_1,j_2,\cdots,j_m}$, since every $\theta_{j_1,j_2,\cdots,j_m}$ is composed of these polynomials. How many $x_1^{k_1}x_2^{k_2}\cdots x_n^{k_n}$ meet the requirements? We can estimate this number by counting the ways in which $|G|$ can be expressed as a sum of $n$ non-negative numbers, which is 
    \begin{equation}
	    \frac{(|G|+n-1)!}{(n-1)!|G|!}.
	\end{equation}
	This formula shows exponentially growing complexity, so it is not economical to reduce the precision of the numerical roots by increasing the complexity.
	
	All $\theta_{j_1,j_2,\cdots,j_i}$ are polynomials of $x_1, x_2, \cdots, x_n$, and the number of terms increases exponentially with the increase of $i$, but if they are regarded as polynomials of $\theta_{j_1,j_2,\cdots,j_{(i-1)}}$, the number of terms remains below a fixed value. Therefore, we should substitute the approximate value of $\theta_{j_1,j_2,\cdots,j_{(i-1)}}$ into $\theta_{j_1,j_2,\cdots,j_i}$ in each loop to get the approximate values of $\theta_{j_1,j_2,\cdots,j_i}$. Thus we can avoid saving $\theta_{j_1,j_2,\cdots,j_i}$ on the computer as polynomials of $x_1,x_2,\cdots,x_n$. One difficulty is that $\theta_{j_1,j_2,\cdots,j_i}$ not only depends on $\theta_{j_1,j_2,\cdots,j_{(i-1)}}$, but also on $\sigma^k_i\theta_{j_1,j_2,\cdots,j_{(i-1)}}$. The calculation for the latter requires the polynomial forms of $\theta_{j_1,j_2,\cdots,j_{(i-1)}}$. In order to solve this problem, we need to calculate $p_1p_2\cdots p_{i}$ values of $\theta_{j_1,j_2,\cdots,j_i}$ under the action of $G/G_i$ in each loop: $\sigma_m^{k_m}\sigma_{m-1}^{k_{m-1}}\cdots\sigma_{i+1}^{k_{i+1}}\theta_{j_1,j_2,\cdots,j_i}$. The number of these values is only $p_1p_2\cdots p_m=|G|$ at most. When all $\sigma_m^{k_m}\sigma_{m-1}^{k_{m-1}}\cdots\sigma_{i}^{k_i}\theta_{j_1,j_2,\cdots,j_{(i-1)}}$ have been calculated, since $\theta_{j_1,j_2,\cdots,j_i}$ is regarded as polynomials of $\sigma_i^k\theta_{j_1,j_2,\cdots,j_{(i-1)}}$, all $\sigma_m^{k_m}\sigma_{m-1}^{k_{m-1}}\cdots\sigma_{i+1}^{k_{i+1}}\theta_{j_1,j_2,\cdots,j_i}$ can be calculated. Therefore, the loop of the algorithm can continue, and the exponential growth of the number of polynomial terms is avoided. After each calculation of the approximate values of $\sigma_m^{k_m}\sigma_{m-1}^{k_{m-1}}\cdots\sigma_{i+1}^{k_{i+1}}\theta_{j_1,j_2,\cdots,j_i}$, we continue to calculate the approximate value of each $(\zeta_i^k,\theta_{j_1,j_2,\cdots,j_{(i-1)}})$, then all the approximations of $\sigma_m^{k_m}\sigma_{m-1}^{k_{m-1}}\cdots\sigma_{i}^{k_i}\theta_{j_1,j_2,\cdots,j_{(i-1)}}$ can be deleted for freeing memory. We save each approximation of $(\zeta_i^k,\theta_{j_1,j_2,\cdots,j_{(i-1)}})$ in order to determine the uncertain phase caused by the $p_{i}$th root extraction in the step from $(\zeta_i^k,\theta_{j_1,j_2,\cdots,j_{(i-1)}})^{p_i}$ to $(\zeta_i^k,\theta_{j_1,j_2,\cdots,j_{(i-1)}})$.
	
	Based on this idea, we can further reduce the quantity of calculation of the algorithm. Since
	\begin{align}
	    \begin{split}
	        (\zeta_i^k,\theta_{j_1j_2,\cdots, j_{(i-1)}})^{p_i} =& \theta_{j_1,j_2,\cdots,j_{(i-1)},0} + \theta_{j_1,j_2,\cdots,j_{(i-1)},1}\zeta_i^k + \theta_{j_1,j_2,\cdots,j_{(i-1)},2}\zeta_i^{2k} +\cdots \\
	        & \quad + \theta_{j_1,j_2,\cdots,j_{(i-1)},(p_i-1)}\zeta_i^{(p_i-1)k},
	    \end{split}
	\end{align}
	similar to Eq.~\eqref{eq:sigmax}, we have
	\begin{align}\label{eq:theta1i}
	    \begin{split}
	        \theta_{j_1,j_2,\cdots,j_{(i-1)},j_i} = \frac{1}{p_i}\Big[&(\zeta_i^0,\theta_{j_1j_2,\cdots, j_{(i-1)}})^{p_i}+\zeta_i^{-j_i}(\zeta_i^1,\theta_{j_1j_2,\cdots, j_{(i-1)}})^{p_i}+\\ &\zeta_i^{-2j_i}(\zeta_i^2,\theta_{j_1j_2,\cdots, j_{(i-1)}})^{p_i}+\cdots+\zeta_i^{-(p_i-1)j_i}(\zeta_i^{p_i-1},\theta_{j_1j_2,\cdots, j_{(i-1)}})^{p_i}\Big].
	    \end{split}
	\end{align}
	With the action of $\sigma_m^{k_m}\sigma_{m-1}^{k_{m-1}}\cdots\sigma_{i+1}^{k_{i+1}}$ on both sides of Eq.~\eqref{eq:theta1i}, we can still get a correct equation. Thus, we can avoid dealing with polynomial expansions and go straight to the numerical result we need.
    
    Based on the above analysis, I summarize the whole algorithm as follows: 
    
    1. Find the numerical roots that satisfy the requirement of Eq.~\eqref{eq:accuate}, and find the approximate value of each primitive $p_i$th root of unity with the same significant digits as the numerical roots. 
    
    2. Take any root $x_1$, and use the approximate value of each root to construct a $p_1\times p_1\times\cdots\times p_m$-dimensional array $\Theta_0$ which satisfies
    \begin{equation}
	    \Theta_0[j_1,j_2,\cdots,j_m] = \sigma_m^{j_m}\sigma_{m-1}^{j_{m-1}}\cdots\sigma_2^{j_2}\sigma_1^{j_1}x_1,
	\end{equation}
	where $0\leq j_i < p_i$. Since the action of the group $G$ on the roots is transitive, $\Theta_0$ contains all the roots. 
	
	3. The cycle of the program begins. In each loop, $i$ takes an integer value from 1 to $m$ in turn, and we calculate the following $p_1\times p_1\times\cdots\times p_m$-dimensional array $L_{i-1}$:
	\begin{equation}\label{eq:Li}
	    L_{i-1}[j_1,\cdots,j_{i-1},k,j_{i+1},\cdots,j_m] = \sum_{j=0}^{p_i-1}\zeta_i^{jk}\Theta_{i-1}[j_1,\cdots,j_{i-1},j,j_{i+1},\cdots,j_m].
	\end{equation}
	The relationship between the array $L_{i-1}$ and the Lagrange resolvents is as follows:
	\begin{equation}
	    \sigma_m^{j_m}\cdots\sigma_{i+1}^{j_{i+1}}(\zeta_i^k,\theta_{j_1j_2,\cdots, j_{(i-1)}}) = L_{i-1}[j_1,\cdots,j_{i-1},k,j_{i+1},\cdots,j_m].
	\end{equation}
	With the array $L_{i-1}$, we can calculate the array $\Theta_i$:
	\begin{align}\label{eq:Thetai}
	    \Theta_i[j_1,\cdots,j_{i-1},j_i,j_{i+1},\cdots,j_m]=\frac{1}{p_i}\sum_{k=0}^{p_i-1}\zeta_i^{-kj_i}(L_{i-1}[j_1,\cdots,j_{i-1},k,j_{i+1},\cdots,j_m])^{p_i}.
	\end{align}
	For this $\Theta_i$, we have
	\begin{equation}
	    \sigma_m^{j_m}\cdots\sigma_{i+1}^{j_{i+1}}\theta_{j_1,j_2,\cdots,j_{(i-1)},j_i} = \Theta_i[j_1,\cdots,j_{i-1},j_i,j_{i+1},\cdots,j_m].
	\end{equation}
	The main calculation amount of the whole algorithm is concentrated in this step, and the main part is the multiplications. Therefore, the complexity of the algorithm can be obtained by counting the number of multiplications in the entire loop. There are $p_i|G|$ multiplications in Eq.~\eqref{eq:Li}. In Eq.~\eqref{eq:Thetai}, calculating the $p_i$th power of all elements of the array $L_{i-1}$ requires $(p_i-1)|G|$ times of multiplication, and then calculating the elements of $\Theta_i$ requires $p_i|G|$ multiplications, so the number of multiplications required by Eq.~\eqref{eq:Thetai} is $(2p_i-1)|G|$. Therefore, the total number of multiplications required for the this step is 
	\begin{equation}
	    |G|\sum_{i=1}^m(3p_i-1)<3|G|^2
	\end{equation}
	On the other hand, the required significant digits of the numerical roots is proportional to |G|, which makes the computation amount of a single multiplication that maintains the accuracy is proportional to $|G|^2$. Therefore, the complexity of the algorithm is less than $A |G|^4$. 
	
	4. After getting the array $\Theta_m$, round the value of each element of $\Theta_m$ to get the exact value of the corresponding $\theta_{j_1,j_2,\cdots,j_m}$. Then the radical formulas for the roots of the equation can be obtained by  working backward step by step.

    \section{One Simple Example: \texorpdfstring{$x^5+20x+32=0$}{x5+20x+32=0}}
    \label{sec:example}
    
    The Galois group $G$ of $f(x)=x^5+20x+32$ is the dihedral group $D_5$~\cite{JENSEN1982347}, and the generators are $(1,2,3,4,5)$ and $(1,4)(2,3)$ when $G$ is represented by permutation. The composition series of $G$ is~\cite{GAP4}
    \begin{equation}
	    G=G_2\rhd G_1=\langle(1,2,3,4,5)\rangle\rhd G_0=\{e\},
	\end{equation}
	The generator of $G_2/G_1$ is $(1,4)(2,3)G_1$, so here we take $\sigma_1=(1,2,3,4,5)$, $\sigma_2=(1 ,4)(2,3)$. The magnitude of the root of $f(x)=0$ is $2.4$, $|G|=2\times 5=10$, we have 
	\begin{equation}
	    \log_{10}(2\times5^{5\times2}\times2^2\times10\times2.4^{10-1})=12.3,
	\end{equation}
	Therefore here we take the approximations of the roots to be accurate to the 13th decimal place. The approximations for these roots are
	\begin{equation}\label{eq:rootsnum}
	\begin{split}
	    \tilde{x}_1 &= -1.3639621650899; \\
	    \tilde{x}_2 &= -1.1078748900075 - 1.7187891044417i; \\
	    \tilde{x}_3 &= 1.7898559725525 + 1.5514288842038i; \\
	    \tilde{x}_4 &= 1.7898559725525 - 1.5514288842038i; \\
	    \tilde{x}_5 &= -1.1078748900075 + 1.7187891044417i.
	\end{split}
	\end{equation}
	I use tilde to differentiate between approximate and exact values. The approximations in Eq.~\eqref{eq:rootsnum} are sorted so that they match the permutation representation of $G$. Since $(2,4,5,3)G(2,4,5,3)^{-1}=G$, there are other possible orderings.
	
	Define
	\begin{equation}
	    \Theta_0=\begin{bmatrix}
	                   \tilde{x}_1,\, \sigma_2\tilde{x}_1 \\
	                   \sigma_1\tilde{x}_1,\, \sigma_2\sigma_1\tilde{x}_1\\
	                   \sigma_1^2\tilde{x}_1,\, \sigma_2\sigma_1^2\tilde{x}_1\\
	                   \sigma_1^3\tilde{x}_1,\, \sigma_2\sigma_1^3\tilde{x}_1 \\
	                   \sigma_1^4\tilde{x}_1,\, \sigma_2\sigma_1^4\tilde{x}_1
	             \end{bmatrix}=
	             \begin{bmatrix}
	                 \tilde{x}_1,\, \tilde{x}_4\\
	                 \tilde{x}_2,\, \tilde{x}_3\\
	                 \tilde{x}_3,\, \tilde{x}_2\\
	                 \tilde{x}_4,\, \tilde{x}_1\\
	                 \tilde{x}_5,\, \tilde{x}_5
	             \end{bmatrix}.
	\end{equation}
	Using Eq.~\eqref{eq:Li} and Eq.~\eqref{eq:Thetai}, the array $\Theta_2$ can be obtained after two loops (retaining 14 significant digits and ignoring the imaginary part): 
	\begin{align*}
	    \Theta_2[0,0] &= 1.4863999240547\times10^{-19} \Rightarrow \theta_{0,0}=0, \\
	    \Theta_2[0,1] &= 1.4863999240547\times10^{-19} \Rightarrow \theta_{0,1}=0, \\
	    \Theta_2[1,0] &= -9999999.9999970 \Rightarrow \theta_{1,0} = -10000000,  \\
	    \Theta_2[1,1] &= 34999999.999995 \Rightarrow \theta_{1,1} = 35000000, \\
	    \Theta_2[2,0] &= 9999999.9999970 \Rightarrow \theta_{2,0} = 10000000,  \\
	    \Theta_2[2,1] &= 14999999.999999 \Rightarrow \theta_{2,1} = 15000000, \\
	    \Theta_2[3,0] &= 9999999.9999970 \Rightarrow \theta_{3,0} = 10000000,  \\
	    \Theta_2[3,1] &= 14999999.999999 \Rightarrow \theta_{3,1} = 15000000, \\
	    \Theta_2[4,0] &= -9999999.9999970 \Rightarrow \theta_{4,0} = -10000000,  \\
	    \Theta_2[4,1] &= 34999999.999995 \Rightarrow \theta_{4,1} = 35000000.
	\end{align*}
	With these results, we can then work back step by step to find the exact expressions for the roots. The arrays $L_0$ and $L_1$ generated during the loops are not listed here. When we encounter the square root or fifth root extraction, we need $L_1$ and $L_0$ to determine the phases.

    \section{Conclusions}
    \label{sec:conclusions} 
    The method presented in this work can systematically solve solvable polynomial equations with rational coefficients of arbitrary degree when we know the corresponding Galois group. In particular, with the help of numerical roots obtained by existing efficient numerical root-finding algorithms, we can effectively reduce the complexity of finding the exact radical formulas of the roots of the solvable polynomial equations. The efficiency of this algorithm benefits from this fact: for the Galois group $G$ of the monic integer coefficient polynomial $f(x)$, any multivariate $G$-invariant polynomial with integer coefficients takes an integer value after substituting the roots of $f(x)$. The whole method does not need to process any polynomials, but only needs to perform basic arithmetic operations. The complexity of the entire algorithm is a polynomial function of $|G|$ ($\leq A|G|^4$). For those solvable polynomials with small Galois groups~\cite{JENSEN1982347,BRUEN1986305,Leonardo2000}, we can even solve them manually with the help of this algorithm and the corresponding numerical roots. However, when the order of the Galois group is large, high demands are placed on the precision of the numerical roots.

    \addcontentsline{toc}{section}{Acknowledgments}
    \acknowledgments
    Thanks for the wonderful life.

    \addcontentsline{toc}{section}{References}
    \bibliographystyle{JHEP}
    \bibliography{bibliography}
    
\end{document}